\documentstyle[11pt]{article}
 \newlength{\standardunitlength}
\newtheorem{cor}{Corollary} \newtheorem{lemma}{Lemma}
\newtheorem{theorem}{Theorem} \newtheorem{prop}{Proposition}
\newenvironment{proof}{\noindent {\sc Proof:}}{$\Box$ \vspace{2 ex}}

\begin{document}

\begin{center} Applications of Symmetric Functions to Cycle and
Increasing Subsequence Structure after Shuffles \end{center}

\begin{center}
Version 1: February 22, 2001
\end{center}

\begin{center}
Referee suggestions implemented on 1/26/02
\end{center}

\begin{center}
By Jason Fulman
\end{center}

\begin{center}
Current affiliation: University of Pittsburgh
\end{center}

\begin{center}
Department of Mathematics
\end{center}

\begin{center}
301 Thackeray Hall
\end{center}

\begin{center}
Pittsburgh, PA 15260
\end{center}

\begin{center}
email:fulman@math.pitt.edu
\end{center}

\begin{center}
Affiliation at time of writing: Stanford University
\end{center}

\newpage \begin{abstract} Using symmetric function theory, we study
the cycle structure and increasing subsequence structure of permutations after iterations of various shuffling methods. We emphasize the role of Cauchy type identities and variations of the Robinson-Schensted-Knuth correspondence. \end{abstract} 

Keywords: Card shuffling, RSK correspondence, cycle index, increasing
subsequence. 

\newpage

\section{Introduction}

	In an unpublished effort to study the way real people shuffle
cards, Gilbert-Shannon-Reeds introduced the following model, called
$k$-riffle shuffling. Given a deck of $n$ cards, one cuts it into $k$
piles with probability of pile sizes $j_1,\cdots,j_k$ given by
$\frac{{n \choose j_1,\cdots,j_k}}{k^n}$. Then cards are dropped from
the packets with probability proportional to the pile size at a given
time (thus if the current pile sizes are $A_1,\cdots,A_k$, the next
card is dropped from pile $i$ with probability
$\frac{A_i}{A_1+\cdots+A_k}$).
	
	The theory of riffle shuffling is relevant to many parts of
mathematics. One area of mathematics influenced by shuffling is Markov
chain theory \cite{D}. For instance Bayer and Diaconis \cite{BD}
proved that $\frac{3}{2}log_2(n)$ $2$-shuffles are necessary and
sufficient to mix up a deck of $n$ cards and observed a cut-off
phenomenon. The paper \cite{H} gives applications of shuffling to
Hochschild homology and the paper \cite{BergeronWolfgang} describes
the relation with explicit versions of the Poincar\'e-Birkhoff-Witt
theorem. Section 3.8 of \cite{SSt} describes GSR shuffles in the
language of Hopf algebras. In recent work, Stanley \cite{Sta} has
related biased riffle shuffles with the Robinson-Schensted-Knuth
correspondence, thereby giving an elementary probabilistic
interpretation of Schur functions and a different approach to some
work of interest to the random matrix community. He recasts many of
the results of \cite{BD} and \cite{F1} using quasisymmetric
functions. Connections of riffle shuffling with dynamical systems
appear in \cite{BD}, \cite{La1}, \cite{La2},
\cite{F4}. Generalizations of the GSR shuffles to other Coxeter groups
appear in \cite{BB},\cite{F2}, \cite{F3}, \cite{F4}, \cite{F5}.

	It is useful to recall one of the most remarkable
properties of GSR $k$-shuffles. Since $k$-shuffles induce a
probability measure on conjugacy classes of $S_n$, they induce a
probability measure on partitions $\lambda$ of $n$. Consider the
factorization of random degree $n$ polynomials over a field $F_q$ into
irreducibles. The degrees of the irreducible factors of a randomly
chosen degree $n$ polynomial also give a random partition of $n$. The
fundamental result of Diaconis-McGrath-Pitman (DMP) \cite{DMP} is that
this measure on partitions of $n$ agrees with the measure induced by
card shuffling when $k=q$. This allowed natural questions on shuffling
to be reduced to known results on factors of polynomials and vice
versa. Lie theoretic formulations, generalizations, and analogs of the
DMP theorem appear in \cite{F2},\cite{F3},\cite{F4}.

	The motivation behind this paper was to understand the DMP
theorem and its cousins in terms of symmetric function theory. (All
notation will follow that of \cite{Mac} and background will appear in
Section \ref{background}). For the DMP theorem itself Stanley
\cite{Sta} gives an argument using ideas from symmetric theory. The
argument in Section \ref{Rif} is different and emphasizes the role of
the RSK correspondence and the Cauchy identity \[ \sum_{\lambda}
s_{\lambda}(x) s_{\lambda}(y) = \sum_{\lambda} \frac{1}{z_{\lambda}}
p_{\lambda}(x) p_{\lambda}(y). \] Here $s_{\lambda}$ and $p_{\lambda}$
denote the Schur functions and power sum symmetric functions
respectively.

 Given Section \ref{Rif}, it was very natural to seek card shuffling interpretations for the Cauchy type identities

\[ \sum_{\lambda} s_{\lambda'}(x) s_{\lambda}(y) = \sum_{\lambda}
\frac{\epsilon_{\lambda}}{z_{\lambda}} p_{\lambda}(x) p_{\lambda}(y) \]

\[ \sum_{\lambda} s_{\lambda}(x) S_{\lambda}(y) =\sum_{\lambda \atop all \ parts \ 
odd} \frac{2^{l(\lambda)}}{z_{\lambda}} p_{\lambda}(x)
p_{\lambda}(y) \]

\[ \sum_{\lambda} s_{\lambda'}(x) S_{\lambda}(y) =\sum_{\lambda \atop all \ parts \ 
odd} \frac{2^{l(\lambda)} \epsilon_{\lambda}}{z_{\lambda}}
p_{\lambda}(x) p_{\lambda}(y) \]

\[ \sum_{\lambda} s_{\lambda}(x) \tilde{s}_{\lambda}(\alpha,\beta,\gamma) =
\sum_{\lambda} \frac{1}{z_{\lambda}} p_{\lambda}(x)
\tilde{p}_{\lambda} (\alpha,\beta,\gamma) \] Here $\lambda'$ denotes
the transpose of a partition and
$\epsilon_{\lambda}=(-1)^{|\lambda|-l(\lambda)}$ where $l(\lambda)$ is
the number of parts of $\lambda$. $S_{\lambda}$ is a symmetric
function studied for instance by Stembridge \cite{Ste} and defined in
Section \ref{Unimodal}. The symmetric function
$\tilde{s}_{\lambda}(\alpha,\beta,\gamma)$ is an extended Schur function to be
discussed in Section \ref{Extended}. (The fourth identity is actually
a generalization of the second identity though it will be helpful to
treat them differently).

	In fact these identities (and probably many identities from symmetric function theory) are related to card shuffling. Section \ref{Bottom} relates the first of these identities to riffle shuffles followed by reversing the order of the cards; the resulting cycle index permits calculations of interest to
real-world shufflers. Section \ref{Unimodal} relates the second of these identities to the cycle structure of affine hyperoctahedral shuffles, which are
generalizations of unimodal permutations; the third identity shows
that dealing from the bottom of the deck has no effect for these
shuffles. This gives a non-Lie theoretic approach to some results in
\cite{F4} and proves a more general assertion. Although there is some
overlap with the preprint \cite{T} for the case of unimodal
permutations, even in that case the treatment here is quite different
and forces into consideration a variation of the RSK
correspondence, which we believe to be new. We should also point out that Gannon \cite{Ga} was the first to solve the problem of counting unimodal permutations by cycle structure, using completely different ideas. (His results are not in the form of a cycle index and it would be interesting to understand the results in this paper by his technique).

	Section \ref{Extended} develops preliminaries related to the case of extended Schur functions. It defines models of card shuffling called $(\vec{\alpha},\vec{\beta},\gamma)$ shuffles (which include the GSR shuffles) and explains how they iterate. This model contains other shuffles of interest such as iterations of the following procedure. Given a deck of $n$ cards, cut the deck into two piles where the sizes are $k,n-k$ with probability $\frac{{n \choose k}}{2^n}$; then shuffle the size $k$ pile thoroughly and riffle it with the remaining cards. This special case was first studied in \cite{DFP} (their work was on convergence rates, not in cycle structure or increasing subsequence structure). Section \ref{Extended} proves that if one applies the usual RSK correspondence to a permutation distributed as a $(\vec{\alpha},\vec{\beta},\gamma)$ shuffle, then the probability of getting any recording tableau of shape $\lambda$ is the extended Schur function $\tilde{s}_{\lambda}(\vec{\alpha},\vec{\beta},\gamma)$. (When $\gamma \neq 0$ this is equivalent to a result of Kerov/Vershik \cite{KV} and Berele/Remmel \cite{BR}. However the case $\gamma \neq 0$ (which arises for the shuffle in this paragraph), is treated incorrectly in \cite{KV} and not at all in \cite{BR}). 

	Section \ref{Shuffling} applies the results of Sections
\ref{Rif} and \ref{Extended} to find formulas for cycle structure
after $(\vec{\alpha},\vec{\beta},\gamma)$ shuffles; for instance it is
proved that after such a shuffle on a deck of size $n$, the expected
number of fixed points is the sum of the first $n$ extended power sum
symmetric functions evaluated at the relevant parameters. An upper
bound on the convergence rate of these shuffles is derived. Section
\ref{Shuffling} closes with a discussion of convolutions of top to
random shuffles, and remarks that for sufficiently large $n$, $5/6
log_2(n)+c$ 2-riffle shuffles bring the longest increasing subsequence
to its limit distribution.

\section{Background} \label{background}

	This section collects the facts from symmetric function theory
which will be needed later. Chapter 1 of \cite{Mac} is a superb
introduction to symmetric functions. We review a few essentials here.

	The power sum symmetric functions $p_{\lambda}$ are an
orthogonal basis of the ring of symmetric functions. Letting
$z_{\lambda}= \prod_i i^{n_i} n_i!$ be the centralizer size of the conjugacy of $S_n$
indexed by the partition $\lambda$ with $n_i$ parts of size $i$, one has that
\[ \langle p_{\lambda},p_{\mu} \rangle = \delta_{\lambda \mu} z_{\lambda}.\] The
descent set of a permutation $w$ is defined as the set of $i$ with $1
\leq i \leq n-1$ such that $w(i) > w(i+1)$; the ascent set is the set
of $i$ with $1 \leq i \leq n-1$ such that $w(i) < w(i+1)$. The descent
set of a standard Young tableau $T$ is the set of $i$ such that $i+1$
is in a lower row of $T$ than $i$. The RSK correspondence (carefully
exposited in \cite{Sa},\cite{SVol2}) associates to a permutation $w$ a pair of
standard Young tableau (its insertion tableau $P(w)$ and its recording
tableau $Q(w)$) and the descent set of $w$ is equal to the descent set
of $Q(w)$. Further the descent set of $w^{-1}$ is equal
to the descent set of $P(w)$, since $Q(w^{-1})=P(w)$. Des$(w)$ and
Asc$(w)$ will denote the descent and ascent set of $w$
respectively. The notation $\lambda \vdash n$ means that $\lambda$ is
a partition of $n$. The symbol $f_{\lambda}$ denotes the number of
standard Young tableau of shape $\lambda$.

	The following result is a simple consequence of work of Gessel
and Reutenauer \cite{GR} and Garsia \cite{Gars}.

\begin{theorem} \label{GRH} Let $\beta_{\lambda}(D)$ be the number of standard
Young tableau of shape $\lambda$ with descent set $D$. Let $N_i(w)$ be
the number of $i$-cycles of a permutation $w$. Then

\begin{enumerate}

\item \[ \sum_{w \in S_n \atop Des(w)=D} \prod_{i \geq 1} x_i^{N_i(w)} =
\langle \sum_{\lambda \vdash n} s_{\lambda}(y) \beta_{\lambda}(D), \prod_{i,j
\geq 1} e^{\frac{x_i^j}{ij} \sum_{d|i} \mu(d) p_{jd}(y)^{i/d}} \rangle \]

\item \[ \sum_{w \in S_n \atop Asc(w)=D} \prod_{i \geq 1} x_i^{N_i(w)} =
\langle \sum_{\lambda \vdash n} s_{\lambda}'(y) \beta_{\lambda}(D), \prod_{i,j
\geq 1} e^{\frac{x_i^j}{ij} \sum_{d|i} \mu(d) p_{jd}(y)^{i/d}} \rangle.\]

\end{enumerate}
\end{theorem} 

\begin{proof} The number of $w$ in $S_n$ with descent set $D$ and
$n_i$ $i$-cycles is the coefficient of $\prod_i x_i^{N_i(w)}$ on the
left hand side of the first equation. Let $\tau$ be the partition with
$n_i$ parts of size $i$ and let Lie$_{\tau}(y)$ be the symmetric function
associated with the corresponding Lie character (for background on Lie characters and relevant symmetric function theory see \cite{R}). By \cite{GR}, the
number of $w$ in $S_n$ with descent set $D$ and $n_i$ $i$-cycles is
equal to the inner product

\[ \langle \sum_{\lambda \vdash n} s_{\lambda}(y) \beta_{\lambda}(D),
Lie_{\tau}(y) \rangle .\] From \cite{Gars} it follows that Lie$_{\tau}(y)$ is
the coefficient of $\prod_i x_i^{n_i}$ in \[ \prod_{i,j \geq 1}
e^{\frac{x_i^j}{ij} \sum_{d|i} \mu(d) p_{jd}(y)^{i/d}}.\] This proves
the first assertion.

	For the second assertion, note that $\beta_{\lambda'}(D) =
\beta_{\lambda}(\{1,\cdots,n-1\}-D)$. This follows from the fact that
if a permutation $w$ has RSK shape $\lambda$ and descent set $D$, then
its reversal has RSK shape $\lambda'$ and ascent set $D$. Thus

\begin{eqnarray*}
&  & \langle \sum_{\lambda \vdash n} s_{\lambda}'(y) \beta_{\lambda}(D), \prod_{i,j
\geq 1} e^{\frac{x_i^j}{ij} \sum_{d|i} \mu(d)
p_{jd}(y)^{i/d}} \rangle\\
& = & \langle \sum_{\lambda \vdash n} s_{\lambda}(y) \beta_{\lambda}(\{1,\cdots,n-1\}-D), \prod_{i,j
\geq 1} e^{\frac{x_i^j}{ij} \sum_{d|i} \mu(d)
p_{jd}(y)^{i/d}} \rangle
\end{eqnarray*} as desired.
\end{proof}

\section{Biased riffle shuffles} \label{Rif}

	We emphasize from the start that the main result in this
subsection is not new: it is equivalent to assertions proved in
\cite{F1} and then in work of Stanley \cite{Sta}. It was first proved
for ordinary riffle shuffles in \cite{DMP}. The value of the current
argument is that it underscores the role of RSK and the Cauchy
identity

\[ \sum_{\lambda} s_{\lambda}(x) s_{\lambda}(y) = \sum_{\lambda}
\frac{1}{z_{\lambda}} p_{\lambda}(x) p_{\lambda}(y) \] (the sums are
over all partitions of all natural numbers). 

	Biased riffle shuffles were introduced in \cite{DFP} and
studied further in \cite{F1}. A biased riffle shuffle with parameters
$\vec{q}=(q_1,q_2,\cdots)$ where $\sum q_i=1$ is defined as
follows. First cut the deck into piles of sizes $k_1,k_2,\cdots$ by
picking the $k$'s according to the distribution \[ {n \choose k_1,k_2,
\cdots} \prod_{i} q_i^{k_i}.\] Now drop cards from the packets one at
a time, according to the rule that at each stage the probability of
dropping from a packet is proportional to the number of cards in that
packet. For instance if there are 2 packets with sizes $3$ and $5$,
then the next card would come from the first packet with probability
$3/8$. It is not hard to see that the probability that a biased riffle
shuffle gives a permutation $w$ depends on $w$ only through
Des$(w^{-1})$. The main case of interest is $q_1=\cdots=q_k=1/k$ all
other $q_i=0$ and corresponds to ordinary riffle shuffles \cite{BD}.

	To determine the cycle structure after a biased riffle shuffle
we could make use of the following result of Stanley \cite{Sta}.

\begin{theorem} \label{Rich} Let $w$ be distributed as a biased riffle
shuffle with parameters $\vec{q}$. Let $T$ be a standard Young tableau
of shape $\lambda$. Then the probability that the RSK algorithm
associates insertion tableau $T$ to $w$ is equal to
$s_{\lambda}(\vec{q})$. \end{theorem} Instead (to simplify later
sections) we will use the following similar result, which we record
for completeness.

\begin{theorem} \label{record1} Let $w$ be distributed as a biased
riffle shuffle with parameters $\vec{q}$. Let $T$ be a standard Young
tableau of shape $\lambda$. Then the probability that the RSK
algorithm associates recording tableau $T$ to $w$ is equal to
$s_{\lambda}(\vec{q})$. \end{theorem}

\begin{proof} Given a length $n$ word $J$ on the symbols $\{1, \cdots,
k\}$, let $a_i$ be the number of occurrences of symbol $i$ in $J$
respectively. Define a permutation $w$ in two line form by putting
$1,\cdots,a_1$ in the positions occupied by the $1$'s of $J$ from left
to right, then putting the next $a_2$ numbers in the positions
occupied by the $2$'s of $J$ from left to right, and so on. For
instance the word \[ 1 \ 3 \ 2 \ 1 \ 2 \ 2 \ 1 \ 3 \ 1 \ 2 \]
corresponds to the permutation \[ 1 \ 9 \ 5 \ 2 \ 6 \ 7 \ 3 \ 10 \ 4 \
8.\] It is easy to see that in general the recording tableau of $w$
under the RSK algorithm is equal to the recording tableau of $J$ under
the RSK algorithm. Arguing as in \cite{BD}, if the entries of the
random word $J$ are chosen independently with probability $q_i$ of
symbol $i$, then the resulting distribution on permutations $w$ is the
same as performing a $\vec{q}$ biased riffle shuffle. As in \cite{KV},
the combinatorial definition of the Schur function immediately implies
that the chance that $J$ has recording tableau $T$ is
$s_{\lambda}(\vec{q})$. \end{proof}

	Lemma \ref{annoying} could be simplified via Theorem
\ref{Rich} but we prefer not to take this path.

\begin{lemma} \label{annoying} Let $\beta_{\lambda}(D)$ be the number
of standard Young tableau of shape $\lambda$ with descent set $D$. If
$\beta_{\lambda}(D) \neq 0$, then the probability that a biased
$\vec{q}$-shuffle produces a permutation $w$ with Des$(w^{-1})=D$ and
RSK shape $\lambda$ is equal to the probability that a biased
$\vec{q}$-shuffle produces a permutation with $(P,Q)$ tableaux
satisfying Des$(P(w))=D$,shape$(Q(w))=\lambda$ divided by
$\beta_{\lambda}(D) f_{\lambda}$. \end{lemma}

\begin{proof} Fix any permutation $w$ such that Des$(w^{-1})=D$ and
such that $w$ has RSK shape $\lambda$ (this is possible if
$\beta_{\lambda}(D) \neq 0$). Let $x$ be the probability of obtaining
$w$ after a biased $\vec{q}$ shuffle. Since all $w$ with
Des$(w^{-1})=D$ are equally likely, $x=y/z$ where $y$ is the
probability that a biased $\vec{q}$ shuffle leads to a permutation
with inverse descent set $D$ and RSK shape $\lambda$, and $z$ is the
number of permutations with inverse descent set $D$ and RSK shape
$\lambda$. Now $y$ is the probability that after a biased $\vec{q}$
shuffle one obtains a permutation $w$ with Des$(P(w))=D$,
shape$(Q(w))=\lambda$. Note that $z$ is simply $\beta_{\lambda}(D)
f_{\lambda}$, since the insertion tableau can be any standard Young
tableau of shape $\lambda$ and descent set $D$, and the recording
tableau can be any standard Young tableau of shape
$\lambda$. \end{proof}

	Now we prove the main result in this subsection.

\begin{theorem} \label{cyc1} Let $E_{n,\vec{q}}$ denote expected value
under the biased riffle shuffle measure with parameters
$\vec{q}$. Let $N_i(w)$ be the number of $i$-cycles of the permutation $w$. Then

\[ \sum_{n \geq 0} u^n E_{n,\vec{q}} (\prod_i x_i^{N_i}) =
\prod_{i,j} e^{\frac{(u^ix_i)^j}{ij} \sum_{d|i} \mu(d)
p_{jd}(\vec{q})^{i/d}}.\]

\end{theorem}

\begin{proof} Let $w$ be a fixed permutation such that
Des$(w^{-1})=D$; then Prob$_{\vec{q}}(D)$ will denote the probability
of obtaining $w$ after a biased riffle shuffle with parameters
$\vec{q}$. Using part 1 of Theorem \ref{GRH} one concludes that the
sought cycle index is

\begin{eqnarray*}
& & \sum_{n \geq 0} u^n E_{n,\vec{q}} (\prod_i x_i^{N_i})\\
& = & \sum_{n \geq 0} u^n \sum_{n_i \geq 0}(\prod_i x_i^{n_i}) \sum_{D
\subseteq \{1,\cdots,n-1\}} Prob_{\vec{q}}(D) |\{w:Des(w)=D,N_i(w)=n_i\}|\\
& = & \sum_{n \geq 0} \sum_{D \subseteq \{1,\cdots,n-1\}} Prob_{\vec{q}}(D) \langle \sum_{\lambda \vdash n} s_{\lambda}(y) \beta_{\lambda}(D), \prod_{i,j
\geq 1} e^{\frac{(u^ix_i)^j}{ij} \sum_{d|i} \mu(d)
p_{jd}(y)^{i/d}} \rangle\\
& = & \sum_{n \geq 0} \langle \sum_{\lambda \vdash n} s_{\lambda}(y)  \sum_{D \subseteq \{1,\cdots,n-1\}} Prob_{\vec{q}}(D)
\beta_{\lambda}(D) ,
\prod_{i,j \geq 1} e^{\frac{(u^ix_i)^j}{ij} \sum_{d|i} \mu(d)
p_{jd}(y)^{i/d}} \rangle
\end{eqnarray*}

	Lemma \ref{annoying} implies that \[ \sum_{D \subseteq
\{1,\cdots,n-1\}} Prob_{\vec{q}}(D) \beta_{\lambda}(D) \] is
$\frac{1}{f_{\lambda}}$ multiplied by the probability that the
recording tableau of a permutation obtained after a biased $\vec{q}$
shuffle has shape $\lambda$. By Theorem \ref{record1}, this latter
probability is $s_{\lambda}(\vec{q}) f_{\lambda}$. Hence the sought
cycle index is simply the inner product

\[ \langle \sum_{\lambda} s_{\lambda}(y) s_{\lambda} (\vec{q}), \prod_{i,j
\geq 1} e^{\frac{(u^ix_i)^j}{ij} \sum_{d|i} \mu(d)
p_{jd}(y)^{i/d}} \rangle.\] Applying the Cauchy identity yields

\[ \langle \sum_{\lambda} \frac{1}{z_{\lambda}} p_{\lambda}(y) p_{\lambda}
(\vec{q}), \prod_{i,j \geq 1} e^{\frac{(u^ix_i)^j}{ij} \sum_{d|i}
\mu(d) p_{jd}(y)^{i/d}} \rangle.\] Since $\langle p_{\lambda},p_{\mu} \rangle = \delta_{\lambda,\mu} z_{\lambda}$ this simplifies to

\[ \prod_{i,j \geq 1} e^{\frac{(u^ix_i)^j}{ij} \sum_{d|i} \mu(d)
p_{jd}(\vec{q})^{i/d}}.\] \end{proof} 

	We remark that for $k$-riffle shuffles the cycle index
simplifies to

\[ \prod_{i \geq 1} (\frac{1}{1-\frac{u^ix_i}{k^i}})^{\frac{1}{i}
\sum_{d|i} \mu(d) k^{i/d}}.\]

\section{Dealing from the Bottom of the Deck} \label{Bottom}

	This section considers cycle structure of a biased riffle
shuffle followed by dealing from the bottom of the deck. This is
equivalent to turning the deck upside down after shuffling. (Persi
Diaconis points out that someone running card guessing experiments
might do this). The results in this section are all new. Results about
subsequence structure are omitted since reversing the order of a
permutation simply transposes its RSK shape.

	Let $\lambda'$ denote the transpose of $\lambda$. Let
$l(\lambda)$ be the number of parts of $\lambda$ and let
$\epsilon_{\lambda}$ denote $(-1)^{|\lambda|-l(\lambda)}$. Whereas the
previous subsection used the Cauchy identity, this subsection uses the
dual Cauchy identity

\[ \sum_{\lambda} s_{\lambda'}(x) s_{\lambda}(y) = \sum_{\lambda}
\frac{\epsilon_{\lambda}}{z_{\lambda}} p_{\lambda}(x) p_{\lambda}(y)
\] (the sums are over all partitions of all natural numbers).

\begin{theorem} \label{cyclereverse} Let $E'_{n,\vec{q}}$ denote
expected value under the biased riffle shuffle measure with parameters
$\vec{q}$ followed by reversing the order of the cards. Then

\[ \sum_{n \geq 0} u^n E'_{n,\vec{q}} (\prod_i x_i^{N_i}) =
\prod_{i,j} e^{\frac{((-u)^ix_i)^j}{ij} \sum_{d|i} \mu(d)
(-p_{jd}(\vec{q}))^{i/d}}.\]

\end{theorem}

\begin{proof} Let $w$ be a fixed permutation such that Asc$(w^{-1})=D$;
then Prob$'_{\vec{q}}(D)$ will denote the probability of obtaining $w$
after a $\vec{q}$ biased riffle shuffle followed by reversing the
order of the cards.

	Using part 2 of Theorem \ref{GRH} one concludes that the
sought cycle index is

\begin{eqnarray*}
& & \sum_{n \geq 0} u^n E'_{n,\vec{q}} (\prod_i x_i^{N_i})\\
& = & \sum_{n \geq 0} u^n \sum_{n_i \geq 0}(\prod_i x_i^{n_i}) \sum_{D
\subseteq \{1,\cdots,n-1\}} Prob'_{\vec{q}}(D)
|\{w:Asc(w)=D,N_i(w)=n_i\}|\\
& = & \sum_{n \geq 0} \sum_{D \subseteq \{1,\cdots,n-1\}}
Prob'_{\vec{q}}(D) \langle \sum_{\lambda \vdash n} s_{\lambda'}(y)
\beta_{\lambda}(D), \prod_{i,j \geq 1} e^{\frac{(u^ix_i)^j}{ij}
\sum_{d|i} \mu(d) p_{jd}(y)^{i/d}} \rangle\\
& = & \sum_{n \geq 0} \langle \sum_{\lambda \vdash n} s_{\lambda'}(y)  \sum_{D \subseteq \{1,\cdots,n-1\}} Prob'_{\vec{q}}(D)
\beta_{\lambda}(D),
\prod_{i,j \geq 1} e^{\frac{(u^ix_i)^j}{ij} \sum_{d|i} \mu(d)
p_{jd}(y)^{i/d}} \rangle. \end{eqnarray*}

	From the proof of Theorem \ref{cyc1},
\[ \sum_{D \subseteq \{1,\cdots,n-1\}} Prob'_{\vec{q}}(D) \beta_{\lambda}(D)  = \sum_{D \subseteq \{1,\cdots,n-1\}} Prob_{\vec{q}}(D) \beta_{\lambda}(D) = s_{\lambda}(\vec{q}).\] Consequently the sought cycle index is simply the inner product

\[ \langle \sum_{\lambda} s_{\lambda'}(y) s_{\lambda} (\vec{q}), \prod_{i,j
\geq 1} e^{\frac{(u^ix_i)^j}{ij} \sum_{d|i} \mu(d) p_{jd}(y)^{i/d}} \rangle.\]
Applying the dual Cauchy identity yields

\[ \langle \sum_{\lambda} \frac{\epsilon_{\lambda}}{z_{\lambda}} p_{\lambda}(y)
p_{\lambda} (\vec{q}), \prod_{i,j \geq 1} e^{\frac{(u^ix_i)^j}{ij}
\sum_{d|i} \mu(d) p_{jd}(y)^{i/d}} \rangle.\] Since $\langle p_{\lambda},p_{\mu} \rangle = \delta_{\lambda,\mu} z_{\lambda}$ this simplifies to

\[ \prod_{i,j \geq 1} e^{\frac{(u^ix_i)^j}{ij} \sum_{d|i} \mu(d)
(-1)^{ji-i/d} p_{jd}(\vec{q})^{i/d}}.\]
\end{proof} 
	
	The case of most interest is $q_1=\cdots=q_k=\frac{1}{k}$ and
all other $q_i=0$. Then the cycle index simplifies to

\[ \prod_{i \geq 1} (\frac{1}{1-\frac{(-u)^ix_i}{k^i}})^{\frac{1}{i} \sum_{d|i}
\mu(d) (-k)^{i/d}}.\] Much information can be gleaned from this cycle index in
analogy with results in \cite{DMP} for ordinary riffle shuffles
(i.e. when one deals from the top of the deck). We record three such
results which are perhaps the the most interesting.

\begin{cor} \label{fixedpoints} The expected number of fixed points
after a $k$-riffle shuffle on n cards followed by reversing the order of the
cards is \[ 1 - \frac{1}{k} + \frac{1}{k^2} \cdots +
\frac{(-1)^{n-1}}{k^{n-1}}.\] \end{cor}

\begin{proof} The generating function for fixed points is given by
setting $x_i=1$ for all $i>1$ in the cycle index. This yields \[
(1+x_1u/k)^k \prod_{i \neq 1 }
(\frac{1}{1-\frac{(-u)^i}{k^i}})^{\frac{1}{i} \sum_{d|i} \mu(d)
(-k)^{i/d}}.\] Multiplying and dividing by $(1+u/k)^k$ gives \[
\frac{(1+x_1u/k)^k}{(1+u/k)^k} \prod_{i}
(\frac{1}{1-\frac{(-u)^i}{k^i}})^{\frac{1}{i} \sum_{d|i} \mu(d)
(-k)^{i/d}}.\] Observe that \[ \frac{1}{1-u} = \prod_{i \geq 1}
(\frac{1}{1-\frac{(-u)^i}{k^i}})^{\frac{1}{i} \sum_{d|i} \mu(d)
(-k)^{i/d}}\] since this is what one obtains by setting all $x_i=1$ in
the cycle index. Hence the generating function for fixed points is \[
\frac{(1+x_1u/k)^k}{(1+u/k)^k (1-u)}.\] Then one differentiates with
respect to $x_1$, sets $x_1=1$, and takes the coefficient of
$u^n$. \end{proof}

	We remark that \cite{DMP} showed that the expected number of
fixed points for $k$-riffle shuffles on an n-card deck is \[ 1 +
\frac{1}{k} + \frac{1}{k^2} \cdots + \frac{1}{k^{n-1}}.\] It is
straightforward to compute higher moments for $k$ shuffles followed by
reversal.

	The next goal is to determine the limit behavior of the
distributions of the short cycles. The answer differs considerably
from the GSR riffle shuffle case, in which only convolutions of
geometric distributions come into play.

	We require a simple lemma.

\begin{lemma} \label{largen} If $f(u)$ has a Taylor series $\sum_{n
\geq 0} a_n u^n$ which converges at $u=1$, then the $n \rightarrow
\infty$ limit of the coefficient of $u^n$ in $\frac{f(u)}{1-u}$ is
$f(1)$. \end{lemma}

\begin{proof} This follows because the coefficient of $u^n$ in
$\frac{f(u)}{1-u}$ is $a_0+\cdots+a_n$. \end{proof}

\begin{cor} \label{limitbehavior} 
\begin{enumerate}

\item Fix $u$ such that $0<u<1$. Choose a random deck size with
probability of getting $n$ equal to $(1-u)u^n$. Let $N_i(w)$ be the
number of $i$-cycles of $w$ distributed as the reversal of a $k$
riffle shuffle. Then the random variables $N_i$ are independent, where
$N_i$ ($i$ odd) is the convolution of $\frac{1}{i} \sum_{d|i} \mu(d)
k^{i/d}$ many binomials with parameter $u^i/(k^i+u^i)$ and $N_i$ ($i$
even) is the convolution of $\frac{1}{i} \sum_{d|i} \mu(d) (-k)^{i/d}$
many geometrics with parameter $u^i/k^i$.

\item Let $N_i(w)$ be the number of $i$-cycles of $w$ distributed as
the reversal of a $k$ riffle shuffle. Then as $n \rightarrow \infty$
the random variables $N_i$ converge in finite dimensional distribution
to independent random variables, where $N_i$ ($i$ odd) becomes the
convolution of $\frac{1}{i} \sum_{d|i} \mu(d) k^{i/d}$ many binomials
with parameter $1/(k^i+1)$ and $N_i$ ($i$ even) becomes the
convolution of $\frac{1}{i} \sum_{d|i} \mu(d) (-k)^{i/d}$ many
geometrics with parameter $1/k^i$.  \end{enumerate} \end{cor}

\begin{proof} As noted after Theorem \ref{cyclereverse}, the cycle
index of a $k$-shuffle followed by reversing the order of the cards is

\[ \prod_{i \geq 1} (\frac{1}{1-\frac{(-u)^ix_i}{k^i}})^{\frac{1}{i}
\sum_{d|i} \mu(d) (-k)^{i/d}}.\] The proof of Corollary
\ref{fixedpoints} gives that

\[ \frac{1}{1-u} = \prod_{i \geq 1}
(\frac{1}{1-\frac{(-u)^i}{k^i}})^{\frac{1}{i} \sum_{d|i} \mu(d)
(-k)^{i/d}}. \] Dividing these equations implies that

\begin{eqnarray*} & & \sum_{n \geq 0} (1-u) u^n E'_{n,1/k,\cdots,1/k}
(\prod_i x_i^{N_i})\\ & = & \prod_{i \ odd}
(\frac{1+u^ix_i/k^i}{1+u^i/k^i})^{1/i \sum_{d|i} \mu(d) k^{i/d}}
\prod_{i \ even} (\frac{1-u^i/k^i}{1-u^ix_i/k^i})^{1/i \sum_{d|i}
\mu(d) (-k)^{i/d}}.  \end{eqnarray*} This proves the first assertion
of the theorem. The second assertion follows from dividing both sides
of this equation by $1-u$ and applying Lemma \ref{largen}. (Note that
if all but finitely many $x_i=1$, only finitely many terms in the
generating function remain. Since $k \geq 2$ the Taylor series
converges at $u=1$ provided that the remaining $x$'s aren't too much
larger than 1). \end{proof}

	Finally we observe (Corollary \ref{largecycles}) that the
distribution of the large cycles is the same as for random
permutations (in contrast to the case of small cycles). One can guess
this heuristically from the generating function since the large $i$
terms of the cycle index converge to those of random permutations. The
same happens for ordinary riffle shuffles (Proposition 5.5 of
\cite{DMP}). The distribution of large cycles in random permutations
has been broadly studied (\cite{VS} and the references therein).

\begin{cor} \label{largecycles} Fix $k$ and let $L_1,\cdots,L_r$ be the lengths of the $r$ longest cycles of $\pi$. Then for $k$ fixed, or growing with $n$ as $n \rightarrow \infty$, \[ |Prob'_{n,1/k,\cdots,1/k}(L_1/n \leq t_1,\cdots,L_r/n \leq t_r) - Prob_{S_n}(L_1/n \leq t_1,\cdots L_r/n \leq t_r)| \rightarrow 0 \] uniformly in $t_1,\cdots,t_r$. (Here $Prob_{S_n}$ denotes the uniform distribution on $S_n$). \end{cor}

\begin{proof} Given the cycle index for k-shuffles followed by a
reversal, this follows from minor modifications of either the
arguments in \cite{Hans} or \cite{ABT}. \end{proof}

\section{Unimodal Permutations and a Variation of the RSK
Correspondence} \label{Unimodal}

	One goal of this section is to understand cycle structure
after shuffling by the following method.

\begin{center}
Generalized Shuffling Method on $C_n$
\end{center}

 	Step 1: Start with a deck of $n$ cards face down. Let $0 \leq
y_1,\cdots,y_k \leq 1$ be such that $\sum y_i=1$. Choose numbers
$j_1,\cdots,j_{2k}$ multinomially with the probability of getting
$j_1,\cdots,j_{2k}$ equal to ${n \choose j_1,\cdots,j_{2k}}
\prod_{i=1}^{k} y_i^{j_{2i-1}+j_{2i}}$. Make $2k$ stacks of cards of
sizes $j_1,\cdots,j_{2k}$ respectively. Flip over the even numbered
stacks.

	Step 2: Drop cards from packets with probability proportional
to packet size at a given time. Equivalently, choose uniformly at
random one of the ${n \choose j_1,\cdots,j_{2k}}$ interleavings of the
packets.

	Cycle structure of this model of shuffling was analyzed for
equal $y$ in \cite{F4}. (Actually there one flipped over the odd
numbered piles, but this has no effect on the cycle index as the
resulting sums in the group algebra are conjugate by the longest
element in $S_n$. By a result of Sch$\ddot{u}$tzenberger exposited as
Theorem A1.2.10 in \cite{SVol2}, conjugation by the longest element
also has no effect on RSK shape). The model was introduced for $k=1$
(and thus $y_1=1$) in \cite{BD}. Let $E^*_{n,\vec{y}}$ be expectation
on $C_n$ after the above shuffling method. Let $N_i(w)$ be the number
of $i$-cycles of $w$ in $C_n$, disregarding signs. It is proved in
\cite{F4} that

\begin{theorem} \label{Brau} \begin{eqnarray*} & & 1+\sum_{n \geq 1}
u^n \sum_{w \in C_n} E^*_{n,\frac{1}{k},\cdots,\frac{1}{k}} (\prod_{i
\geq 1} x_i^{N_i(w)})\\ & = & \prod_{m \geq 1} (\frac{1+x_m
u^m/(2k)^m}{1-x_m u^m/(2k)^m})^{\frac{1}{2m} \sum_{d|m \atop d \ odd}
\mu(d) (2k)^{\frac{m}{d}}}.  \end{eqnarray*} \end{theorem} As the
paper \cite{F4} did not discuss asymptotics of long cycles, before
proceeding we note the following corollary, whose proof method is the
same as that of Corollary \ref{largecycles}.

\begin{cor} Fix $k$ and let $L_1,\cdots,L_r$ be the lengths of the $r$ longest cycles of $\pi$. Then for $k$ fixed, or growing with $n$ as $n \rightarrow \infty$, \[ |Prob^*_{n,1/k,\cdots,1/k}(L_1/n \leq t_1,\cdots,L_r/n \leq t_r) - Prob_{S_n}(L_1/n \leq t_1,\cdots L_r/n \leq t_r)| \rightarrow 0 \] uniformly in $t_1,\cdots,t_r$. (Here $Prob_{S_n}$ denotes the uniform distribution on $S_n$). \end{cor}

	A generalization of Theorem \ref{Brau} will be proved later in this
section. To this end, we require the following variation of the RSK
correspondence.

{\bf Variation of the RSK Correspondence:} Order the set of numbers
$\{\pm 1, \cdots , \pm k\}$ by \[ 1 < -1 < 2 < -2 \cdots < k < -k.\]
Given a word on these symbols, run the RSK algorithm as usual, with
the amendments that a symbol $i$ can't bump another $i$ if $i$ is positive,
but must bump another $i$ if $i$ is negative. (This guarantees that positive
numbers appear at most once in each column and that negative numbers
appear at most once in each row).

	For example the word \[ 1 \ -1 \ 2 \ -2 \ 1 \ 1 \ -1 \ 1 \ 2 \
2 \ -1 \ 2 \ -2 \] has insertion tableau $P$ and recording tableau $Q$
respectively equal to

\[ \begin{array}{c c c c c c c c c c c}
             & & & 1 & 1 & 1 & 1 & -1  & 2 &2  &-2 \\
            &&  &  -1 & 2 & 2 &  &&& &\\
             &&&    -1 & -2 &  &  & &&&  \end{array}  \]

\[ \begin{array}{c c c c c c c c c c c}
             & & & 1 & 2 & 3 & 4 & 9  & 10 & 12  &13 \\
            &&   &  5 & 6 & 7 & && & &\\
             &&&    8 & 11 &  &  & &&&  \end{array}  \] 

	The proof of Theorem \ref{vary} runs along the same lines as
the proof of the RSK correspondence as presented in
\cite{Sa}. Hence we omit the details.

\begin{theorem} \label{vary} Order the set of numbers $\{\pm 1, \cdots
, \pm k\}$ by \[ 1 < -1 < 2 < -2 \cdots < k < -k.\] Then the above
variation on the RSK Correspondence is a bijection between length $n$
words on the symbols $\{\pm 1, \cdots , \pm k\}$ and pairs $(P,Q)$
where

\begin{enumerate}

\item $P$ is a tableau on the symbols $\{\pm 1, \cdots, \pm k\}$
satisfying $P(a,b) \leq P(a+1,b)$, $P(a,b) \leq P(a,b+1)$ for all
$a,b$ where $P(a,b)$ denotes the entry in the $a$th row and $b$th
column of $P$.

\item If $i$ is positive then it appears at most once in each column
of $P$ and if $i$ is negative then it appears at most once in each row
of $P$.

\item $Q$ is a standard Young tableau on the symbols $\{1,\cdots,n\}$.

\item $P$ and $Q$ have the same shape.

\end{enumerate}
\end{theorem} 	

	The next result relates the shuffling model of this section
with the above variation of the RSK correspondence. For its statement,
$S_{\lambda}$ will denote the symmetric functions studied in
\cite{Ste} (a special case of the extended Schur functions in
\cite{KV}). One definition of the $S_{\lambda}$ is as the determinant

\[ S_{\lambda}(y) = \det(q_{\lambda_i-i+j}) \] where $q_{-r}=0$ for
$r>0$ and for $r \geq 0$, $q_r$ is defined by setting \[ \sum_{n \geq
0} q_nt^n = \prod_{i \geq 1} \frac{1+y_it}{1-y_it}.\] We remark that
Theorem \ref{probinterp} gives a simple probabilistic interpretation
to $S_{\lambda}$, different from the interpretation in \cite{KV}.

\begin{theorem} \label{probinterp} Let $w$ be distributed as a shuffle
of this section with parameters $y_1,\cdots,y_k$ after forgetting
about signs. Let $Q$ be a standard Young tableau of shape
$\lambda$. Then the probability that the usual RSK correspondence
associates recording tableau $Q$ to $w$ is equal to $\frac{1}{2^n}
S_{\lambda}(y_1,\cdots,y_k)$. Consequently the probability that $w$
has RSK shape $\lambda$ is equal to $\frac{f_{\lambda}}{2^n}
S_{\lambda}(y_1,\cdots,y_k)$. \end{theorem}

\begin{proof} Given a length $n$ word $J$ on the symbols $\{\pm 1,
\cdots, \pm k\}$, let $a_i,b_i$ be the number of occurrences of the
symbol $i,-i$ in $J$ respectively. Define a permutation $w$ in two line
form by putting $1,\cdots,a_1$ in the positions occupied by the $1$'s
of $J$ from left to right, then putting the next $b_1$ numbers
(arranged in decreasing order) in the positions occupied by the $-1$'s
of $J$ from left to right, then the next $a_2$ numbers (arranged in
increasing order) in the positions occupied by the $2$'s of $J$ from
left to right, etc. For instance the word \[ 1 \ -1 \ 2 \ -2 \ 1 \ 1 \
-1 \ 1 \ 2 \ 2 \ -1 \ 2 \ -2 \] corresponds to the permutation \[ 1 \
7 \ 8 \ 13 \ 2 \ 3 \ 6 \ 4 \ 9 \ 10 \ 5 \ 11 \ 12.\] If the word
entries are chosen independently with $\pm i$ having probability
$\frac{y_i}{2}$, the resulting distribution on permutations is the same
as performing a $\vec{y}$ shuffle of this section and forgetting about
signs.

	It is easy to see that the recording tableau of $w$ under the
RSK algorithm is equal to the recording tableau of $J$ under our
variant of the RSK algorithm. Let $\gamma_i(P)$ be the number of
occurrences of symbol $i$ in a tableau $P$.  By Theorem \ref{vary},
the probability that $J$ has recording tableau $Q$ under our variant
of RSK is equal to

\[ \frac{1}{2^n} \sum_{P} \prod_{i \geq 1}
y_i^{\gamma_i(P)+\gamma_{-i}(P)} \] where $P$ has shape $\lambda$ and
satisfies conditions 1,2 in Theorem \ref{vary}. Theorem 9.2b of
\cite{Ste} shows that this sum is equal to $\frac{1}{2^n}
S_{\lambda}(y_1,\cdots,y_k)$. \end{proof}

	As mentioned in the introduction, Theorem \ref{probinterp}
is relevant to random matrix theory. This is
because the first row in the RSK shape of a random permutation $w$ is
equal to the length of the longest increasing subsequence of $w$ and
has asymptotically the same distribution as the largest eigenvalue of
a random GUE matrix \cite{BDJ}. Studying longest increasing
subsequences of $w$ distributed as a GSR $k$-riffle shuffle amounts to
studying the longest weakly increasing subsequences in random length
$n$ words on $k$ symbols, which has also been of interest to random
matrix theorists \cite{Sta,TW}. What Theorem \ref{probinterp} tells
us is that studying longest increasing subsequences of $w$ distributed
as unsigned type $C$ shuffles amounts to studying weakly increasing
subsequences in random length $n$ words on the symbols $\{ \pm 1,
\cdots, \pm k\}$, where $1<-1<\cdots<k<-k$ and the subsequence is not
allowed to contain a given negative symbol $i$ more than once. For $k$ fixed
and random length $n$ words on the symbols $\{ \pm 1, \cdots,
\pm k\}$, roughly half the symbols will be positive, and the negative
symbols can in total affect the length of the longest weakly
increasing subsequence by at most $k$. For example, one obtains the following corollary from the analogous results in \cite{J} and \cite{TW} for weakly increasing subsequences in random words.

\begin{cor} For $k$ fixed, the RSK shape after an unsigned $C_n$
shuffle with $y_1=\cdots=y_k=\frac{1}{k}$ has at most $k$ rows and $k$
columns. For large $n$ the expected value of any of the $k$ rows or
columns is asymptotic to $\frac{n}{2k}$. \end{cor} We hope in future
work to study the fluctuations around this limit shape, and to examine
the case when both $n,k$ are large.
	
	Theorem \ref{cycleunimodal} determines the generating function
for cycle structure after performing the generalized shuffling method
on $C_n$ with parameters $y_1,\cdots,y_k$ and forgetting about signs.

\begin{theorem} \label{cycleunimodal} Let $E^*_{n,\vec{y}}$ denote
expected value under the generalized shuffling method on $C_n$ with
parameters $y_1,\cdots,y_k$ after forgetting signs. As usual, let $N_i(\pi)$ be the number of cycles of length $i$ of the permutation $\pi$. Then

\[ \sum_{n \geq 0} u^n E^*_{n,\vec{y}} (\prod_i x_i^{N_i}) = \prod_{i
\geq 1} \prod_{j \ odd} e^{\frac{(u^ix_i/2^i)^j}{ij} \sum_{d|i \atop d
\ odd} \mu(d) (2 p_{jd}(y))^{i/d}}.\] Furthermore, reversing the order
of the cards has no effect on the cycle index. \end{theorem}

\begin{proof} Let $w$ be a fixed permutation such that Des$(w^{-1})=D$ and
let Prob$^*_{\vec{y}}(D)$ be the probability of obtaining $w$ after a
$\vec{y}$ unsigned type $C$ shuffle.

	Using part 1 of Theorem \ref{GRH} and the fact that the
probability of $w$ depends only on $w$ through Des$(w^{-1})$, it
follows that the sought cycle index is

\begin{eqnarray*}
& & \sum_{n \geq 0} u^n E^*_{n,\vec{y}} (\prod_i x_i^{N_i})\\
& = & \sum_{n \geq 0} u^n \sum_{n_i \geq 0}(\prod_i x_i^{n_i}) \sum_{D
\subseteq \{1,\cdots,n-1\}} Prob^*_{\vec{y}}(D)
|\{w:Des(w)=D,N_i(w)=n_i\}|\\
& = & \sum_{n \geq 0} \sum_{D \subseteq \{1,\cdots,n-1\}}
Prob^*_{\vec{y}}(D) \langle \sum_{\lambda \vdash n} s_{\lambda}(z)
\beta_{\lambda}(D), \prod_{i,j \geq 1} e^{\frac{(u^ix_i)^j}{ij}
\sum_{d|i} \mu(d) p_{jd}(z)^{i/d}} \rangle\\
& = & \sum_{n \geq 0} \langle \sum_{\lambda \vdash n} s_{\lambda}(z)  \sum_{D \subseteq \{1,\cdots,n-1\}} Prob^*_{\vec{y}}(D) \beta_{\lambda}(D) ,
\prod_{i,j \geq 1} e^{\frac{(u^ix_i)^j}{ij} \sum_{d|i} \mu(d)
p_{jd}(z)^{i/d}} \rangle. \end{eqnarray*}

	Arguing as in Theorem \ref{cyc1} shows that \[ \sum_{D
\subseteq \{1,\cdots,n-1\}} Prob^*_{\vec{y}}(D) \beta_{\lambda}(D) =
\frac{1}{2^n} S_{\lambda}(y).\] Thus the sought cycle index is simply the inner product

\[ \langle \sum_{\lambda} s_{\lambda}(z) S_{\lambda} (y), \prod_{i,j
\geq 1} e^{\frac{(u^ix_i/2^i)^j}{ij} \sum_{d|i} \mu(d) p_{jd}(z)^{i/d}} \rangle.\]
Applying the third identity in the introduction (due to Stembridge \cite{Ste}) yields

\[ \langle \sum_{\lambda \atop all \ parts \ odd}
\frac{2^{l_{\lambda}}}{z_{\lambda}} p_{\lambda}(z) p_{\lambda} (y),
\prod_{i,j \geq 1} e^{\frac{(u^ix_i/2^i)^j}{ij} \sum_{d|i} \mu(d)
p_{jd}(z)^{i/d}} \rangle.\] Since $\langle p_{\lambda},p_{\mu} \rangle =
\delta_{\lambda,\mu} z_{\lambda}$, this simplifies as desired to

\[ \prod_{i} \prod_{j \ odd} e^{\frac{(u^ix_i/2^i)^j}{ij} \sum_{d|i
\atop d \ odd} \mu(d) (2 p_{jd}(y))^{i/d}}.\]

	For the second assertion, Theorem \ref{cyclereverse} shows that the cycle
index after reversing the card order at the end is given by

\[ \prod_{i} \prod_{j \ odd} e^{\frac{((-u)^ix_i/2^i)^j}{ij} \sum_{d|i
\atop d \ odd} \mu(d) (-2 p_{jd}(y))^{i/d}}.\] It is easy to see that
the $-$ signs all drop out.
\end{proof}

	We remark that in the case of greatest interest
($y_1=\cdots=y_k=\frac{1}{k}$, all other $y_i=0$), one recovers
Theorem \ref{Brau}.

	A unimodal permutation $w$ on the symbols $\{1,\cdots,n\}$ is
defined by requiring that there is some $i$ with $1 \leq i \leq n$
such that the following two properties hold:
 \begin{enumerate} \item If $a<b\leq i$, then $w(a)<w(b)$.
  \item If $i \leq a<b$, then $w(a)>w(b)$.
  \end{enumerate} Thus $i$ is where the maximum is achieved, and the
permutations $12\cdots n$ and $nn-1\cdots 1$ are counted as
unimodal. For each fixed $i$ there are ${n-1 \choose i-1}$ unimodal
permutations with maximum $i$, hence a total of $2^{n-1}$ such
permutations. As noted in \cite{Ga}, unimodal permutations are those which avoid
the patterns $213$ and $312$.

	Unimodal permutations are the shuffles of this section in the
case $k=1$ after forgetting about signs; hence Theorem \ref{Brau}
(from \cite{F4}) gives a cycle index for unimodal permutations. The
paper \cite{T}, which appeared in between \cite{F4} and this paper,
obtained a count of unimodal permutations by cycle structure and
position of their maximum, denoted by $max(w)$. We prove an equation
equivalent to Thibon's result \cite{T}. The proof uses the notation
that $m_i(\lambda)$ is the number of parts of $\lambda$ of size $i$.

\begin{theorem} \label{thib} Let $N_i(w)$ be the number of $i$-cycles of a permutation $w$.
\[ 1 + \sum_{n \geq 1} u^n (1+t) \sum_{w \ unimodal} t^{max(w)-1}
\prod_i x_i^{N_i(w)} = \prod_{i,j} e^{\frac{(x_i u^i)^j}{ij}
\sum_{d|i} \mu(d)(t^{jd}-(-1)^{jd})^{i/d}} .\]
\end{theorem}

\begin{proof} A permutation on $n$ symbols is unimodal with maximum at
position $k$ if and only if it has descent set
$k,k+1,\cdots,n-1$. Hence Theorem \ref{GRH} implies that

\begin{eqnarray*}
& &  1 + \sum_{n \geq 1} u^n (1+t) \sum_{w \ unimodal} t^{max(w)-1}
\prod_i x_i^{N_i(w)}\\
& = & \langle 1+(1+t) \sum_{a,b  \geq 0} s_{(a+1,1^b)}(z) t^a u^{a+b+1},  \prod_{i,j
\geq 1} e^{\frac{x_i^j}{ij} \sum_{d|i} \mu(d) p_{jd}(z)^{i/d}} \rangle.
\end{eqnarray*} This can be further simplified using Macdonald's identity (page 49 of
\cite{Mac})
\[ 1+(t+u) \sum_{a,b \geq 0} s_{(a+1,1^b)}(z) t^a u^b = \prod_{i \geq 1}
\frac{1+uz_i}{1-tz_i} \] with $t$ replaced by $tu$ to yield

\begin{eqnarray*}
& & \langle \prod_{i \geq 1} \frac{1+uz_i}{1-tuz_i},  \prod_{i,j
\geq 1} e^{\frac{x_i^j}{ij} \sum_{d|i} \mu(d) p_{jd}(z)^{i/d}} \rangle\\
& = & \langle e^{\sum_{i \geq 1} u^ip_i(z) (t^i-(-1)^i)/i},  \prod_{i,j
\geq 1} e^{\frac{x_i^j}{ij} \sum_{d|i} \mu(d) p_{jd}(z)^{i/d}} \rangle\\
& = & \langle \sum_{\lambda} \frac{p_{\lambda}(z) u^{|\lambda|} \prod_i (t^i-(-1)^i)^{m_i(\lambda)}}{z_{\lambda}},  \prod_{i,j
\geq 1} e^{\frac{x_i^j}{ij} \sum_{d|i} \mu(d) p_{jd}(z)^{i/d}} \rangle\\
& = &  \prod_{i,j} e^{\frac{(x_i u^i)^j}{ij}
\sum_{d|i} \mu(d)(t^{jd}-(-1)^{jd})^{i/d}}.
\end{eqnarray*} Note that we have used the identity

\[ \prod_{i \geq 1} \frac{1}{1-uz_i} = e^{\sum_{i \geq 1} p_i(z)u^i/i}.\]
\end{proof}

\section{Extended Schur functions} \label{Extended}

	The extended complete symmetric functions
$\tilde{h}_k(\alpha,\beta,\gamma)$ are defined by the generating
function \[ \sum_{k=0}^{\infty} \tilde{h}_k(\alpha,\beta,\gamma) z^k =
e^{\gamma z} \prod_{i \geq 1} \frac{1+\beta_i z}{1-\alpha_i z}.\] For
$\lambda=(\lambda_1,\cdots,\lambda_n)$, the extended Schur functions
are defined by \[ \tilde{s}_{\lambda} = det(\tilde{h}_{\lambda_i-i+j})_{i,j=1}^n .\] The extended Schur functions give the characters of the infinite symmetric group and are usefully reviewed in \cite{O}. Observe that $\tilde{s}_{\lambda}$ is obtained from taking the expression for $s_{\lambda}$ as a polynomial in the $h_k$ and replacing $h_k$ by $\tilde{h}_k$. Defining a homomorphism $\Phi$ on symmetric functions by $\Phi(h_k)=\tilde{h}_k$, one sees that any identity for ordinary symmetric functions gives a corresponding identity for extended symmetric functions. That is how one derives the Cauchy identity \[ \sum_{\lambda} s_{\lambda}(x) \tilde{s}_{\lambda}(\alpha,\beta,\gamma) = \sum_{\lambda}
\frac{1}{z_{\lambda}} p_{\lambda}(x) \tilde{p}_{\lambda} (\alpha,\beta,\gamma) \] for extended Schur functions from the usual Cauchy identity (e.g. Example 3.23 of \cite{Mac} for the case $\gamma \neq 0$).

	Since probabilities must be positive, one motivation for interpreting extended Schur functions probabilistically is the following positivity result. \begin{theorem} \label{Thoma's} (\cite{E}) Let $G(z) =
\sum_{k=0}^{\infty} g_k z^k$ be such that $g_0=1$ and all $g_k \geq
0$. Then
\[ det(g_{\lambda_i-i+j})_{i,j=1}^n \geq 0 \] for all partitions
$\lambda$ if and only if
\[ G(z) = e^{\gamma z} \prod_{i \geq 1} \frac{1+\beta_i z}{1-\alpha_i
z} \] where $\gamma \geq 0$ and $\sum \beta_i, \sum \alpha_i$ are
convergent series of positive numbers. \end{theorem} 

	Next we define $(\vec{\alpha},\vec{\beta},\gamma)$ shuffles. We suppose that $\gamma+\sum \alpha_i+\sum \beta_i=1$ and that
$\gamma \geq 0$, $\alpha_i,\beta_i \geq 0$ for all $i$. Using these
parameters, we define a random permutation on $n$ symbols as
follows. First, create a word of length $n$ by choosing letters $n$
times independently according to the rule that one picks $i>0$ with
probability $\alpha_i$, $i<0$ with probability $\beta_i$, and $i=0$
with probability $\gamma$. We use the usual ordering $\cdots < -1 < 0
< 1 < \cdots$ on the integers. Starting with the smallest negative
symbol which appears in the word, let $m$ be the number of times it
appears. Then write $\{1,2,\cdots,m\}$ under its appearances in {\it
decreasing} order from left to write. If the next negative symbol
appears $k$ times write $\{m+1,\cdots,m+k\}$ under its appearances,
again in decreasing order from left to write. After finishing with the
negative symbols, proceed to the $0$'s. Letting $r$ be the number of
$0$'s, choose a random permutation of the relevant $r$ consecutive
integers and write it under the $0$'s. Finally, move to the positive
symbols. Supposing that the smallest positive symbol appears $s$
times, write the relevant $s$ consecutive integers under its
appearances in {\it increasing} order from left to right.

 	The best way to understand this procedure is through an
example. Given the string \[ -2 \ 0 \ 1 \ 0 \ 0 \ 2 \ -1 \ -2 \ -1 \ 1
\] one obtains each of the six permutations

\[ 2 \ 5 \ 8 \ 6 \ 7 \ 10 \ 4 \ 1 \ 3 \ 9 \] \[ 2 \ 5 \ 8 \ 7 \ 6 \ 10
\ 4 \ 1 \ 3 \ 9 \] \[ 2 \ 6 \ 8 \ 5 \ 7 \ 10 \ 4 \ 1 \ 3 \ 9 \] \[ 2 \
6 \ 8 \ 7 \ 5 \ 10 \ 4 \ 1 \ 3 \ 9 \] \[ 2 \ 7 \ 8 \ 5 \ 6 \ 10 \ 4 \
1 \ 3 \ 9 \] \[ 2 \ 7 \ 8 \ 6 \ 5 \ 10 \ 4 \ 1 \ 3 \ 9 \] with
probability $1/6$. In all cases the $1,2$ correspond to the $-2$'s,
the $3,4$ correspond to the $-1$'s, the $8,9$ correspond to the $1$'s
and the $10$ corresponds to the $2$. The symbols $5,6,7$ correspond to
the $0$'s and there are six possible permutations of these symbols. We
call this probability measure on permutations a
$(\vec{\alpha},\vec{\beta},\gamma)$ shuffle.

	The following elementary result (generalizing results in
\cite{BD} and \cite{DFP}) gives physical descriptions of these
shuffles and explains how they convolve. The proof method follows that
of \cite{BD}.

\begin{prop} \label{describe}
\begin{enumerate}

\item A $(\vec{\alpha},\vec{\beta},\gamma)$ shuffle is equivalent to
the following procedure. Cut the $n$ card deck into piles with sizes
$X_i$ indexed by the integers, where the probability of having
$X_i=x_i$ for all $i$ is equal to \[ \frac{n!}{\prod_{i=-
\infty}^{\infty} x_i!}  \gamma^{x_0} \prod_{i >0} \alpha_i^{x_i}
\prod_{i<0} \beta_i^{x_i}. \] The top cards go to the non-empty pile
with smallest index, the next batch of cards goes to the pile with
second smallest index, and so on. Then mix the pile indexed by $0$
until it is a random permutation, and turn upside down all of the
piles with negative indices. Finally, riffle the piles together as in
the first paragraph of the introduction and look at the underlying
permutation (i.e. ignore the fact that some cards are upside down).

\item The inverse of a $(\vec{\alpha},\vec{\beta},\gamma)$ shuffle is
equivalent to the following procedure. Randomly label each card of the
deck, picking label $0$ with probability $\gamma$, label $i>0$ with
probability $\alpha_i$ and label $i<0$ with probability
$\beta_i$. Deal cards into piles indexed by the labels, where cards
with negative or zero label are dealt face down and cards with
positive label are dealt face up. Then mix the pile labeled $0$ so
that it is a random permutation and turn all of the face up piles face
down. Finally pick up the piles by keeping piles with smaller labels
on top.

\item Performing a $(\vec{\alpha},\vec{\beta},\gamma)$ shuffle $k$
times is the same as performing the following shuffle. One cuts into
piles with labels given by $k$-tuples of integers $(z_1,\cdots,z_k)$
ordered according to the following rule:

\begin{enumerate}

\item $(z_1,\cdots,z_k) < (z_1',\cdots,z_k')$ if $z_1<z_1'$.

\item $(z_1,\cdots,z_k) < (z_1',\cdots,z_k')$ if $z_1=z_1' \geq 0$ and
$(z_2,\cdots,z_k) < (z_2',\cdots,z_k')$.

\item $(z_1,\cdots,z_k) < (z_1',\cdots,z_k')$ if $z_1=z_1' < 0$ and
 $(z_2,\cdots,z_k) > (z_2',\cdots,z_k')$.

\end{enumerate} The pile is assigned probability equal to the product
of the probabilities of the symbols in the $k$ tuple. Then the shuffle
proceeds as in part 1, where negative piles (piles where the product
of the coordinates of the $k$ tuple are negative) are turned upside
down and piles with some coordinate equal to 0 are perfectly mixed
before the piles are all riffled together.  \end{enumerate} \end{prop}

{\bf Examples} As an example of Proposition \ref{describe}, consider
an $(\alpha_1,\alpha_2; \beta_1,\beta_2;\gamma)$ shuffle with
$n=11$. For part 1, it may turn out that $X_{-2}=2$, $X_{-1}=1$,
$X_0=3$, $X_1=2$, and $X_2=3$. Then the deck is cut into piles
$\{1,2\}$, $\{3\}$, $\{4,5,6\}$, $\{7,8\}$, $\{9,10,11\}$. The first
two piles are turned upside down and the third pile is completely
randomized, which might yield piles $\{2,1\}$, $\{3\}$, $\{5,4,6\}$,
$\{7,8\}$, $\{9,10,11\}$. Then these piles are riffled together as in
the GSR shuffle. This might yield the permutation

\[ 5 \ 2 \ 7 \ 4 \ 8 \ 9 \ 10 \ 3 \ 1 \ 11 \ 6 .\] The inverse
description (part 2) would amount to labeling cards 2,9 with $-2$,
card 8 with $-1$, cards 1,4,11 with $0$, card 3,5 with $1$, and cards
6,7,10 with $2$, and then mixing the $0$ pile as $4,1,11$. Note that
this leads to the permutation (inverse to the previous permutation)

\[ 9 \ 2 \ 8 \ 4 \ 1 \ 11 \ 3 \ 5 \ 6 \ 7 \ 10 .\]

	As an example of part 3, note that doing a
$(\alpha_1;\beta_1;0)$ shuffle twice does not give a
$(\vec{\alpha},\vec{\beta},\gamma)$ shuffle, but rather gives a
shuffle with 4 piles in the order $(-1,1),(-1,-1),(1,-1),(1,1)$ where
pile 1 has probability $\beta_1 \alpha_1$, pile 2 has probability
$\beta_1 \beta_1$, pile 3 has probability $\alpha_1 \beta_1$ and pile
4 has probability $\alpha_1 \alpha_1$. Piles 1 and 3 are turned upside
down before the riffling takes place. From Section \ref{Unimodal} of this paper one can still analyze the cycle structure and RSK shape of these shuffles
even though they aren't $(\vec{\alpha},\vec{\beta},\gamma)$
shuffles. (Actually Section \ref{Unimodal} of this paper looked at shuffles conjugate to these shuffles by the longest element; this clearly has no effect
on the cycle index and has no effect on the RSK shape by a result of
Sch$\ddot{u}$tzenberger exposited as Theorem A1.2.10 in \cite{SVol2}).

	As another example of part 3, note that a shuffle with
parameters $(\alpha_1;0;\gamma)$ repeated twice gives a shuffle with 4
piles in the order $(0,0),(0,1),(1,0),(1,1)$ where the first 3 piles
are completely mixed before all piles are riffled together. This is
clearly the same as a $(\alpha_1^2;0;1-\alpha_1^2)$ shuffle, agreeing
with Lemma 2.1 of \cite{DFP}.

	Berele and Remmel \cite{BR} and independently Kerov and
Vershik \cite{KV} consider the following analog of the RSK
Correspondence (different from the variation in Section \ref{Unimodal} as the BRKV version uses the standard ordering on the integers). Given a word on the symbols $\{\pm 1,\pm 2,\cdots\}$
one runs the RSK correspondence with the amendments that a negative
symbol is required to bump itself, but that a positive symbol can't
bump itself. For example the word \[ 1 \ -1 \ 2 \ -2 \ 1 \ 1 \ -2 \]
has insertion tableau $P$ and recording tableau $Q$ respectively equal
to

\[ \begin{array}{c c c c c c c c c c c}
             & & & -2 & 1 & 1 &  &   &  &  & \\
            &&  &  -2 & 2 &  &  &&& &\\
             &&&    -1 &  &  &  & &&&\\
             &&& 1 &&&&&&& \end{array}  \]

\[ \begin{array}{c c c c c c c c c c c}
             & & & 1 & 3 & 6 &  &   &  &   & \\
            &&   &  2 & 5 &  & && & &\\
             &&&    4 &  &  &  & &&&\\
              &&& 7 &&&&&&& \end{array}  \]

\begin{theorem} \label{Ber1} (\cite{BR},\cite{KV}) The above variation
on the Robinson-Schensted-Knuth correspondence gives a bijection
between words of length $n$ from the alphabet of integers with the
symbol $i$ appearing $n_i$ times and pairs $(P,Q)$ where

\begin{enumerate}
\item The symbol $i$ occurs $n_i$ times in $P$.
\item The entries of $P$ are weakly increasing in rows and columns.
\item Each positive symbol occurs at most once in each column of $P$
and each negative symbol occurs at most once in each row of $P$.
\item $Q$ is a standard Young tableau on the symbols $\{1,\cdots,n\}$.
\end{enumerate} Furthermore, \[ \tilde{s}_{\lambda}(\vec{\alpha},\vec{\beta},0) = \sum_{P \atop shape(P)=\lambda} \prod_{i > 0} \alpha_i^{n_i(P)} \prod_{i<0} \beta_i^{n_i(P)}.\] \end{theorem}

	Theorem \ref{probinter} and Corollary \ref{c1} connect card
shuffling to the extended Schur functions. When $\alpha=0$, this
result is essentially in \cite{BR} and \cite{KV}. The paper
\cite{KV} states a version of Theorem \ref{Ber1} in which there is
also a parameter $\gamma$ (their Proposition 3), but it is incorrect
for $\gamma \neq 0$ as the following counterexample shows. Setting all
parameters other than $\alpha_1=\alpha$ and $\gamma=1-\alpha$ equal to
0, it follows from the definitions that the extended Schur function
$\tilde{s}_2$ is equal to $\frac{\alpha^2+1}{2}$. But if Proposition 3
of \cite{KV} were correct, it would also equal
$\alpha^2+(1-\alpha)\alpha = \alpha$ since the two words giving a
Young tableau with 1 row of length 2 are $11$ and $01$. In fact as the
$2$ in the denominator of $\frac{\alpha^2+1}{2}$ shows, one can't
interpret the extended Schur functions with $\gamma \neq 0$ in terms
of RSK and words on a finite number of symbols. This accounts for the extra
randomization step (choosing a random permutation for the symbols corresponding the 0's) in our definition of $(\vec{\alpha},\vec{\beta},\gamma)$ shuffles.

	Theorem \ref{probinter} give a probabilistic interpretation of $\tilde{s}_{\lambda}$ for all values of $\gamma$.

\begin{theorem} \label{probinter} Let $\pi$ be distributed as a
permutation under a $(\vec{\alpha},\vec{\beta},\gamma)$ shuffle. Let $Q$
be any standard Young tableaux of shape $\lambda$. Then the
probability that $\pi$ has Robinson-Schensted-Knuth recording tableau
equal to $Q$ is
$\tilde{s}_{\lambda}(\vec{\alpha},\vec{\beta},\gamma)$. \end{theorem}

\begin{proof} First suppose that $\gamma = 0$. As indicated earlier in
this section, each length $n$ word $w$ on the symbols $\{\pm 1, \pm 2,
\cdots \}$ defines a permutation $\pi$. From this construction, it is
easy to see that the recording tableau of $w$ under the BRKV variation
of the RSK algorithm is equal to the recording tableau of $\pi$ under
the RSK algorithm. Thus it is enough to prove that the probability
that the word $w$ has BRKV recording tableau $Q$ is
$\tilde{s}_{\lambda}(\vec{\alpha},\vec{\beta},0)$. This is immediate from
Theorem \ref{Ber1}.
	
	Now the case $\gamma \neq 0$ can be handled by introducing $m$
extra symbols between $0$ and $1$--call them
$1/(m+1),2/(m+1),\cdots,m/(m+1)$ and choosing each with probability
$\gamma / m$. Thus the random word is on $\{\pm 1,\pm 2, \cdots\}$ and
these extra symbols. Each word defines exactly one permutation--the
symbols $1/(m+1),2/(m+1),\cdots,m/(m+1)$ are treated as positive. By
the previous paragraph, the probability of obtaining recording tableau
$Q$ is equal to $\tilde{s}_{\lambda}(\vec{\alpha},\vec{\beta})$ where the
associated $\tilde{h}_{k}$ are defined by \[ \sum_{k=0}^{\infty}
\tilde{h}_k(\alpha,\beta) z^k = (\frac{1}{1-\gamma z /m})^m \prod_{i
\geq 1} \frac{1+\beta_i z}{1-\alpha_i z}.\] As $m \rightarrow \infty$,
this distribution on permutations converges to that of a
$(\vec{\alpha},\vec{\beta},\gamma)$ shuffle, and the generating
function of the $\tilde{h}_k$ converges to \[ \sum_{k=0}^{\infty}
\tilde{h}_k(\alpha,\beta,\gamma) z^k = e^{\gamma z} \prod_{i \geq 1}
\frac{1+\beta_i z}{1-\alpha_i z}.\] \end{proof}

\begin{cor} \label{c1} Let $f_{\lambda}$ be the number of standard
Young tableau of shape $\lambda$. Let $\pi$ be distributed as a
permutation under a $(\vec{\alpha},\vec{\beta},\gamma)$ shuffle. Then
the probability that $\pi$ has Robinson-Schensted-Knuth shape $\lambda$
is equal to $f_{\lambda}
\tilde{s}_{\lambda}(\vec{\alpha},\vec{\beta},\gamma)$. \end{cor}

\section{Convergence Rates and Cycle Index of $(\alpha,\beta,\gamma)$ shuffles} \label{Shuffling}

	First we derive an upper bound on the convergence rate of
$(\vec{\alpha},\vec{\beta},\gamma)$ shuffles to randomness using
strong uniform times as in \cite{DFP}. The
separation distance between a probability $P(\pi)$ and the uniform
distribution $U(\pi)$ is defined as
$max_{\pi}(1-\frac{Q(\pi)}{U(\pi)})$ and gives an upper bound on total
variation distance. Examples of the upper bound of Theorem
\ref{mybound} are considered later.

\begin{theorem} \label{mybound} The separation distance between $k$
applications of a $(\vec{\alpha},\vec{\beta},\gamma)$ shuffle and
uniform is at most \[ {n \choose 2} \left[ \sum_i (\alpha_i)^2+ \sum_i
(\beta_i)^2 \right]^k .\] Thus $k=2log_{\frac{1}{ \sum_i (\alpha_i)^2+
\sum_i (\beta_i)^2 }}n$ steps suffice to get close to the uniform
distribution.  \end{theorem}

\begin{proof} For each $k$, let $A^k$ be a random $n \times k$ matrix
formed by letting each entry equal $i>0$ with probability $\alpha_i$,
$i<0$ with probability $\beta_i$, and $i=0$ with probability
$\gamma$. Let $T$ be the first time that all rows of $A^k$ containing
no zeros are distinct; from the inverse description of
$(\vec{\alpha},\vec{\beta},\gamma)$ shuffles this is a strong uniform
time in the sense of Sections 4B-4D of Diaconis $\cite{Diac}$, since
if all cards are cut in piles of size one the permutation resulting
after riffling them together is random. The separation distance after
$k$ applications of a $(\vec{\alpha},\vec{\beta},\gamma)$ shuffle is
upper bounded by the probability that $T>k$ \cite{AD}. Let $V_{ij}$ be
the event that rows $i$ and $j$ of $A^k$ are the same and contain no
zeros. The probability that $V_{ij}$ occurs is $\left[ \sum_i
(\alpha_i)^2+ \sum_i (\beta_i)^2 \right]^k$. The result follows
because

\begin{eqnarray*}
Prob(T>k) & = & Prob (\cup_{1 \leq i < j \leq n}) V_{ij}\\
& \leq & \sum_{1 \leq i < j \leq n} Prob(V_{ij})\\
& = & {n \choose 2} \left[ \sum_i (\alpha_i)^2+ \sum_i
(\beta_i)^2 \right]^k
\end{eqnarray*}
\end{proof}

	Taking logarithms of the defining identity for $\tilde{h}_k$, one sees that \[ \tilde{p}_1(\vec{\alpha},\vec{\beta},\gamma) =
\sum_i \alpha_i + \sum_i \beta_i + \gamma =1 \] and (for $n \geq 2$)
\[ \tilde{p}_n(\vec{\alpha},\vec{\beta},\gamma) = \sum_i (\alpha_i)^n
+ (-1)^{n+1} \sum_i (\beta_i)^n. \]

	Theorem \ref{cindex} gives a cycle index after
$(\vec{\alpha},\vec{\beta},\gamma)$ shuffles.

\begin{theorem} \label{cindex}

\begin{enumerate}

\item Let $E_{n,(\vec{\alpha},\vec{\beta},\gamma)}$ denote expected
value after a $(\vec{\alpha},\vec{\beta},\gamma)$ shuffle of an $n$
card deck. Let $N_i(\pi)$ be the number of $i$-cycles of a permutation
$\pi$. Then

\[ \sum_{n \geq 0} u^n E_{n,(\vec{\alpha},\vec{\beta},\gamma)}
(\prod_i x_i^{N_i}) = \prod_{i,j} e^{\frac{(u^ix_i)^j}{ij} \sum_{d|i}
\mu(d) \tilde{p}_{jd}(\vec{\alpha},\vec{\beta},\gamma)^{i/d}}.\]

\item Let $E'_{n,(\vec{\alpha},\vec{\beta},\gamma)}$ denote expected
value after a $(\vec{\alpha},\vec{\beta},\gamma)$ shuffle of an $n$
card deck followed by reversing the order of the cards. Then

\[ \sum_{n \geq 0} u^n E'_{n,(\vec{\alpha},\vec{\beta},\gamma)}
(\prod_i x_i^{N_i}) = \sum_{n \geq 0} u^n
E_{n,(\vec{\beta},\vec{\alpha},\gamma)} (\prod_i x_i^{N_i}).\]

\end{enumerate}
\end{theorem}

\begin{proof} Given the results of Section \ref{Extended}, the proof of the first part runs along exactly the same
lines as in the proof of Theorem \ref{cyc1}. The second assertion
follows from the observation that a
$(\vec{\alpha},\vec{\beta},\gamma)$ shuffle followed by reversing the
order of the cards is conjugate (by the longest length element in the
symmetric group) to a $(\vec{\beta},\vec{\alpha},\gamma)$
shuffle. Alternatively, arguing as in the proof of Theorem
\ref{cyclereverse}, one sees that the effect of
reversing the cards on the cycle index of a
$(\vec{\alpha},\vec{\beta},\gamma)$ shuffle is to get

\[ \prod_{i,j} e^{\frac{((-u)^ix_i)^j}{ij} \sum_{d|i} \mu(d)
(-\tilde{p}_{jd}(\vec{\alpha},\vec{\beta},\gamma))^{i/d}}.\] \end{proof}

{\bf Example 1} As a first application of Theorem \ref{cindex}, we
derive an expression for the expected number of fixed points,
generalizing the expression in \cite{DMP}. To get the generating
function for fixed points, one sets $x_2=x_3=\cdots=1$ in the cycle
index. Using the same trick as in \cite{DMP}, the generating function
simplifies to \[ \frac{1}{1-u} \frac{e^{ux \gamma}}{e^{ux}} \prod_{i
\geq 1} \frac{1-u \alpha_i}{1-ux \alpha_i} \frac{1+ux \beta_i}{1+u
\beta_i}.\] Taking the derivative with respect to $x$ and the
coefficient of $u^n$, one sees that the expected number of fixed
points is \[ \gamma + \sum_{j=1}^n [\sum_i (\alpha_i)^j + (-1)^{j+1}
(\beta_i)^j] .\] This is exactly the sum of the first $n$ extended
power sum functions at the parameters
$(\vec{\alpha},\vec{\beta},\gamma)$.

{\bf Example 2} We suppose that $\vec{\beta}=\vec{0}$ and that
$\alpha_1=\cdots=\alpha_q=\frac{1-\gamma}{q}$. Then the cycle index
simplifies to

\[ \prod_{i \geq 1} \left(\frac{1}{1-x_i
(\frac{u(1-\gamma)}{q})^i}\right)^{\frac{1}{i} \sum_{d|i} \mu(d)
q^{i/d}} \prod_{i \geq 1} e^{\frac{u^ix_i (1-(1-\gamma)^i)}{i}}.\] Of
particular interest is the further specialization $q=1$. Then the
cycle index becomes

\[ \frac{1}{1-x_1u(1-\gamma)} \prod_{i \geq 1} e^{\frac{u^ix_i
(1-(1-\gamma)^i)}{i}}.\]

	Recall that a $(1/2,0,1/2)$ shuffle takes a binomial(n,1/2)
number of cards (a binomial(n,1/2) random variable is equal to $k$ with probability ${n \choose k} / 2^n$), thoroughly mixes them, and then riffles them with the remaining cards. Example 3 on page 140 of \cite{DFP} proves (in
slightly different notation) that the iteration of $k$ (1/2,0,1/2)
shuffles is the same as a $((1/2)^k,0,1-(1/2)^k)$ shuffle. They
conclude (in agreement with Theorem \ref{mybound}) that a
$(1/2,0,1/2)$ shuffle takes $log_2(n)$ steps to be mixed, as compared
to $\frac{3}{2} log_2(n)$ for ordinary riffle shuffles. They also
establish a cut-off phenomenon. From the computation of Example 1 one
sees that the expected number of fixed points also drops and that the
mean mixes twice as fast.

	As another example, consider a $(1-1/n,0,1/n)$
shuffle. Heuristically this is like top to random and \cite{DFP}
proves that the convergence rate is the same ($nlog(n)$ steps), which
agrees with Theorem \ref{mybound}. From page 139 of \cite{DFP},
performing a $(1-1/n,0,1/n)$ shuffle $k$ times is the same as
performing a single $((1-1/n)^k,0,1-(1-1/n)^k)$ shuffle. Example 1
gives a formula for the expected number of fixed points. See Example 4 for more discussion of iterations of top to random shuffles.

	Next we consider the asymptotics of cycle structure. As usual, $\mu$ denotes the Moebius function of elementary number theory. Note that considerable simplifications take place when $q=1$ (the interesting case) because $\sum_{d|i} \mu(d)$ is $1$ if $i=1$ and is $0$ otherwise. We omit the details of the proof as they are the same as for the corresponding results in Section \ref{Bottom}.

\begin{cor} \label{limitbehavior2} Suppose that $\vec{\beta}
=\vec{0}$ and $\alpha_1=\cdots=\alpha_q=\frac{1-\gamma}{q}$.
\begin{enumerate} 

\item Fix $u$ such that $0<u<1$. Choose a random deck size with
probability of getting $n$ equal to $(1-u)u^n$. Let $N_i(\pi)$ be the
number of $i$-cycles of $\pi$ distributed as a
$(\vec{\alpha},\vec{\beta},\gamma)$. Then the random variables $N_i$ are independent,
where $N_i$ is the convolution of a Poisson$((u^i(1-(1-\gamma)^i))/i)$
with $\frac{1}{i} \sum_{d|i} \mu(d) q^{i/d}$ many geometrics with
parameter $(\frac{u(1-\gamma)}{q})^i$.

\item Let $N_i(\pi)$ be the number of $i$-cycles of $\pi$ distributed
as a $(\vec{\alpha},\vec{\beta},\gamma)$ shuffle. Then as
$n \rightarrow \infty$ the random variables $N_i$ are independent,
where $N_i$ is the convolution of a Poisson$((1-(1-\gamma)^i)/i)$ with
$\frac{1}{i} \sum_{d|i} \mu(d) q^{i/d}$ many geometrics with parameter
$(\frac{1-\gamma}{q})^i$.

\item Fix $k$ and let $L_1,\cdots,L_r$ be the lengths of the $r$ longest cycles of $\pi$. Then for $k$ fixed, or growing with $n$ as $n \rightarrow \infty$, \[ |Prob'_{n,\vec{\alpha},\vec{\beta},\gamma}(L_1/n \leq t_1,\cdots,L_r/n \leq t_r) - Prob_{S_n}(L_1/n \leq t_1,\cdots L_r/n \leq t_r)| \rightarrow 0 \] uniformly in $t_1,\cdots,t_r$. (Here $Prob_{S_n}$ denotes the uniform distribution on $S_n$).
\end{enumerate}
\end{cor}

{\bf Example 3} Consider the case when
$\alpha_1=\cdots=\alpha_q=\beta_1=\cdots=b_q=\frac{1}{2q}$ and all
other parameters are $0$. Theorems \ref{probinter} and \ref{cindex}
imply that the distribution on RSK shape and cycle index is the same
as for the shuffles in Section \ref{Unimodal}, though we do not see
a simple reason why this should be so. 

{\bf Example 4} Another generalization of riffle shuffles are random walks coming from real hyperplane arrangements \cite{BHR}. The most interesting such shuffles are those where the weights on faces of the Coxeter complex are invariant under the action of the symmetric group. It is straightforward to see that such shuffles are mixtures of what can be called $\mu$ shuffles, where $\mu$ is a composition of $n$. For a $\mu$ shuffle, one breaks the decks into piles of sizes $\mu_1,\mu_2,\cdots$ and then chooses uniformly at random one of the ${n \choose \mu_1,\mu_2, \cdots}$ possible interleavings. In what follows we also let $\mu$ denote the partition of $n$ given by ordering the parts of the composition by decreasing size. 

	For example the top to random shuffle is a $(1,n-1)$
shuffle. Let $P(j,k,n)$ be the probability that when $k$ balls are
dropped at random into $n$ boxes, there are $j$ occupied cells (thus by inclusion exclusion $P(j,k,n)=\sum_{r=j}^n (-1)^{r-j} {n \choose r} {r \choose j}
(1-r/n)^k$). A result of \cite{DFP} is that the iteration of $k$ top
to random shuffles is equivalent to a mixture of $(n-j,1^j)$ shuffles,
where $(n-j,1^j)$ is chosen with probability $P(j,k,n)$. Theorem
\ref{toptorandom} will give an expression for the increasing
subsequence structure after this process. For this a lemma is
required. In its statement we use notation in \cite{Mac} that
$K_{\lambda \mu}$ is a Kostka number (the number of semistandard Young
tableau of shape $\lambda$ where $i$ appears $\mu_i$ times), and
$\lambda/\mu$ denotes a tableau of skew shape $\lambda/\mu$.

\begin{lemma} \label{top} Let $T$ be a standard Young tableau of shape
$\lambda$. The probability that a $\mu$ shuffle has recording tableau
$T$ is equal to $\frac{K_{\lambda \mu}}{{n \choose
\mu_1,\mu_2,\cdots}}$. \end{lemma}

\begin{proof} A $\mu$ shuffle corresponds to choosing at random a word where $i$ appears $\mu_i$ times, and each word has probability ${n \choose \mu_1,\mu_2,\cdots}$. It is easy to see that the RSK recording tableau of the word and the corresponding permutation obtained after the shuffle are identical. Now the number of words of length $n$ where $i$ appears $\mu_i$ times and with recording tableau $T$ is equal to $K_{\lambda \mu}$, since such words biject with the possible insertion tableau which have shape $\lambda$ and weight $\mu$. \end{proof} 

\begin{theorem} \label{toptorandom} Let $f_{\lambda/\mu}$ denote the number of standard tableau of shape $\lambda/\mu$. Then the chance that the RSK shape after $k$ top to random shuffles is $\lambda$ is equal to \[ \frac{f_{\lambda}^2}{n!} \sum_{a=1}^n P(a,k,n) (n-a)! \frac{f_{\lambda/(n-a)}}{f_{\lambda}}.\] \end{theorem} 

\begin{proof} From Lemma \ref{top} and the description of iterations of top to random shuffles as mixtures of $\mu$ shuffles, it follows that the sought probability is \[ \frac{f_{\lambda}}{n!} \sum_{a=1}^n P(a,k,n) K_{\lambda,(n-a,1^a)} (n-a)!.\] Finally observe the $K_{\lambda,(n-a,1^a)}=f_{\lambda/(n-a)}$, since the $n-a$ ones must appear in the first row and what remains is a standard Young tableau. \end{proof}

	Note that in Theorem \ref{toptorandom},
$\frac{f_{\lambda}^2}{n!}$ corresponds to Plancherel measure and the
rest is a correction term (going to 1 as $k \rightarrow \infty$ and
$n$ is fixed). It would be interesting to determine (both for $n$
large and $n$ small) how many iterations of top to random are
necessary for the length of the longest increasing subsequence to be
close to that of a random permutation.

	For comparison, one has the following result for ordinary
2-riffle shuffles. The result is an easy Corollary of equation 1.27 of
\cite{J}, together with the fact that k 2-riffle shuffles is the same
as one $2^k$ riffle shuffle \cite{BD}. Note that the result is for
sufficiently large $n$ and says nothing for $n$ small.

\begin{cor} Let $L_n$ denote the longest increasing subsequence of a
random element of $S_n$ and let $L^(2^k)_n$ denote the longest
increasing subsequence of an element of $S_n$ after k 2-riffle
shuffles. Then \[ lim_{n \rightarrow \infty}
Prob. (\frac{L_n-2n^{1/2}}{n^{1/6}} \leq t) = F(t) \] and \[ lim_{n
\rightarrow \infty} Prob. (\frac{L^{2^k}_n-2n^{1/2}}{n^{1/6}} \leq t)
= F(t-e^{-c})\] where $2^k=\lfloor e^c n^{5/6} \rfloor$ and $F(t)$ is
the Tracy-Widom distribution. Thus for sufficiently large $n$, $5/6
log_2(n)+c$ 2-riffle shuffles are necessary and suffice for the
longest increasing subsequence to be that of a random permutation
. \end{cor}

\section*{Acknowledgements} This research was supported by an NSF
Postdoctoral Fellowship. The author thanks Persi Diaconis and a referee for helpful remarks.

\end{document}